\documentclass[final]{siamart250211}

\usepackage{graphicx} % Required for inserting images
\usepackage{algorithm}
\usepackage{algpseudocode}
\usepackage{amsfonts}
\usepackage{graphicx}
\usepackage{epstopdf}
\ifpdf
  \DeclareGraphicsExtensions{.eps,.pdf,.png,.jpg}
\else
  \DeclareGraphicsExtensions{.eps}
\fi

\newsiamremark{remark}{Remark}
\newsiamremark{hypothesis}{Hypothesis}
\crefname{hypothesis}{Hypothesis}{Hypotheses}
\newsiamthm{claim}{Claim}
\newsiamremark{fact}{Fact}
\crefname{fact}{Fact}{Facts}

\headers{Greedy Thiele interpolation on a continuum}{T. A. Driscoll and Y. Zhou}

\title{Greedy Thiele continued-fraction approximation on continuum domains in the complex plane}

\author{Tobin A. Driscoll\thanks{Department of Mathematical Sciences, University of Delaware, Newark, DE 
  (\email{driscoll@udel.edu}).}
\and Yuxing Zhou\thanks{Department of Mathematical Sciences, University of Delaware, Newark, DE
  (\email{zhouyx@udel.edu}).}
}

\usepackage{amsopn}

\hypersetup{
  pdftitle={Greedy Thiele continued-fraction approximation on continuum domains in the complex plane},
  pdfauthor={T. A. Driscoll and Y. Zhou}
}
\newcommand{\K}{\mathop{\vcenter{\hbox{\scalebox{1.8}{\upshape K}}}}}
\renewcommand{\O}{\mathcal{O}}

\begin{document}

\maketitle

\begin{abstract}
We describe an adaptive greedy algorithm for Thiele continued-fraction (TCF) approximation of a function defined on a continuum domain in the complex plane. The algorithm iteratively selects interpolation nodes from an adaptively refined set of sample points on the domain boundary. We also present new algorithms for evaluating Thiele continued fractions and their accessory weights using only a single floating-point division. Numerical experiments comparing the greedy TCF method with the AAA algorithm on several challenging functions defined on the interval $[-1,1]$ and on the unit circle show that continuum TCF is consistently faster than AAA, by factors ranging from 8 to 40.
\end{abstract}



\section{Introduction}

The use of rational functions, i.e., ratios of polynomials, for approximation is an old subject. But interest in their use as viable numerical objects began to grow rapidly in 2018, with the introduction of the AAA algorithm~\cite{NakatsukasaAAAAlgorithm2018}. This method, based on the barycentric representation of rational functions, allows the construction of rational approximations to a function $f$ defined on a set of sample points in the complex plane. The AAA algorithm adaptively selects a subset of $n$ of the sample points as interpolation nodes and chooses $n$ weights to minimize a linearized residual at the remaining sample points. The AAA algorithm has been extended in various ways~\cite{HochmanFastAAAFast2017,TrefethenComputationZolotarev2025}, and many applications of it are described in the survey~\cite{NakatsukasaApplicationsAAA2026}.

The AAA algorithm requires an SVD at each iteration, and building an interpolant on $n$ nodes (in AAA terminology, support points) takes $\O(n^4)$ work. In practice, the absolute running times are usually well under a second, but the high asymptotic cost is notable. It can be reduced by adopting updating~\cite{HochmanFastAAAFast2017} or random sketching~\cite{ParkFastRandomized2023} techniques, but these are not yet in widespread use. 

Salazar Celis proposed~\cite{SalazarCelisNumericalContinued2024} an alternative greedy iteration for rational approximation based on Thiele's continued fraction (TCF) representation~\cite{ThieleInterpolationsrechnung1909}. This method also adaptively selects nodes from a set of sample points, but it requires no linear algebra and an $\O(n^3)$ cost to build an interpolant on $n$ nodes. The method appears to be numerically stable in practice, although no general proof of this is known.

Both AAA and the greedy TCF method were initially described in terms of fixed sample point sets in the complex plane. In some situations, however, ab initio discretization of the domain boundary is complicated by the presence of unknown nearby singularities in the function to be approximated. Driscoll et al.~\cite{DriscollAAARational2024} described an adaptive refinement strategy for AAA that allows it to work on a continuum domain by adaptively refining the sample set. 

In this paper, we propose a similar continuum strategy for the greedy TCF method. We also point out how to evaluate the TCF representation---and, by extension, the inverse differences analogous to barycentric weights---with a single floating-point division rather than $\O(n)$ of them, resulting in practical speed-ups of 1.6x in the real case and over 6x in the complex case, at least in the Julia programming environment. We also point out that the evaluation is easily augmented to compute derivatives of the interpolant and its residues at simple poles. 

In \autoref{sec:thiele}, we introduce the TCF representation and Salazar Celis' greedy algorithm. In \autoref{sec:continuum}, we describe the new adaptive refinement strategy for continuum domains. In \autoref{sec:newalgs}, we present the new algorithms for evaluation, derivatives, and residues. In \autoref{sec:experiments}, we present numerical experiments comparing the greedy TCF method with AAA on several challenging functions defined on the interval $[-1,1]$ and on the unit circle, demonstrating that TCF is consistently faster than AAA, by factors ranging from 8 to 40. Finally, in \autoref{sec:conclusion}, we summarize our findings.

\section{Thiele's continued fraction representation for rational interpolation}
\label{sec:thiele}

Thiele~\cite{ThieleInterpolationsrechnung1909} introduced a continued fraction representation for rational interpolation in 1909. The classical presentation in~\cite{Milne-ThompsonCalculusFinite1933} is more accessible to English readers. More recent descriptions found in~\cite{Graves-MorrisPracticalReliable1980} and~\cite{GutknechtContinuedFractions1989} generalized the process considerably to cover degeneracies and permit output of approximations with more general degree sequences. Additional context and details can be found in~\cite{JonesContinuedFractions1988} and~\cite{CuytNonlinearMethods1988}, for instance.

Given distinct real or complex nodes $z_1,\ldots, z_n$ at which to take values $y_1,\ldots , y_n$, we define
\begin{equation}
    % r_n(z) = w_1 + \cfrac{z - z_1}{w_2 + \cfrac{z - z_2}{w_3 + \cfrac{z - z_3}{w_4 + \cdots + \cfrac{z-z_{n-1}}{w_n}}}}
    r_n(z) = w_1 + \cfrac{z - z_1}{w_2 + \cfrac{z - z_2}{w_3 + \cfrac{\ddots}{w_{n-1} + \cfrac{z-z_{n-1}}{w_n}}}} = w_1 + \K_{j=1}^{n-1} \frac{z - z_j}{w_{j+1}},
    \label{eq:thiele}
\end{equation}
where the weights $w_k$ are numbers implicitly defined via the interpolation conditions. The continued fraction can also be expressed as $r_n=\psi_1$ resulting from the iteration
\begin{equation}
    \begin{split}
        \psi_n(z) &= w_n, \\
        \psi_{n-k}(z) &= w_{n-k} + \frac{z - z_{n-k}}{\psi_{n-k+1}(z)}, \quad k=1,\ldots,n-1.
    \end{split}
    \label{eq:thiele-eval}
\end{equation}
Requiring that $\psi_1(z_n)=y_n$, we can reverse the iteration to compute
\begin{equation}
    \begin{split}
        % \psi_{n-k+1}(z) &= w_{n-k} + \frac{z_n - z_{n-k}}{\psi_{n-k}(z_n) - w_{n-k}}, \quad k=n-1,n-2,\ldots,1, \\
        \psi_{j+1}(z_n) &= w_{j} + \frac{z_n - z_{j}}{\psi_{j}(z_n) - w_{j}}, \quad j=1,\ldots,n-1, \\
        w_n &= \psi_n(z_n).
    \end{split}
    \label{eq:thiele-weight}
\end{equation}
In fact, the weight computation is itself expressible as the continued fraction
\begin{equation}
    % w_k = 0 + \frac{z_{k} - z_{k-1}}{-w_{k-1}} \cfplus \cdots \cfplus \frac{z_k - z_2}{-w_2} \cfplus \frac{z_k-z_{1}}{-w_1} \cfplus \frac{y_k}{1}, \quad k=1,\dots,n.
    w_n = \K_{j=1}^{n} \frac{z_n - z_{n-j}}{-w_{n-j}}, 
    \label{eq:thiele-weight-cf}
\end{equation}
in which we have made the additional terminal definitions $z_n-z_0=y_n$ and $w_0=-1$. 

Utilizing~\eqref{eq:thiele-weight}, we can start from $(z_1,y_1)$, to find $w_1=y_1$, then iteratively add successive node-value pairs to compute the weights $w_2,\ldots,w_n$ in $\O(n^2)$ operations. If no indeterminate values are encountered in its construction, $r_n(z)$ is a rational function that interpolates the given data. In the absence of cancelling factors in numerator and denominator, $r_n$ is of rational type $(m,m)$ if $n=2m+1$ and type $(m,m-1)$ if $n=2m$. Once the weights have been computed, the evaluation of $r_n(z)$ at a point $z$ can be done in $\O(n)$ operations via, for example,~\eqref{eq:thiele-eval}.

There are two ways the iteration can break down. First, if some $w_n=\infty$, then it must be that $r_{n-1}(z_n)=y_n$---i.e., the most recent interpolant already achieves the desired new value---and the iteration halts. Second, $z_n$ may be an \textit{unattainable} point at which no rational function of the prescribed degree can be constructed. The full characterization of nonexistence is technical (see, e.g.,~\cite{GutknechtContinuedFractions1989}); in practice, such situations are uncommon. Generically, the TCF iteration can be continued in such cases until interpolation is fully restored.

A more troublesome issue is the numerical stability of the computation, which depends greatly on the locations and ordering of the nodes, as well as the values to be interpolated. Cuyt et al.~\cite{CuytInstabilityModification1988} noted that instabilities in the weight computations need not contaminate the evaluation of $r_n$ to the same degree but recommended minimizing approximation degree to avoid instability. Graves–Morris~\cite{Graves-MorrisPracticalReliable1980} showed that a pivoting-style strategy of selecting the smallest-magnitude weight candidate among remaining nodes at each step affords a tractable backward error bound, asserting that the exponential growth of that bound was likely to be a gross overestimate. This assertion was much later proved to be true in the case of Markov functions on real intervals by Beckermann et al.~\cite{BeckermannRationalApproximation2022}. We are not aware of additional numerical or theoretical evidence supporting Graves–Morris' stability assertion. 

%Jones and Thron~\cite{jonesNumericalStability1974} recommended the tail-ordered evaluation $r_n(z)=v_1$ resulting from the iteration
%\begin{equation}
%    \begin{aligned}
%    v_n &= w_n, \\
%    v_k &= w_k + \frac{z-z_k}{v_{k+1}}, \; \text{for } k=n-1,\ldots,1.
%    \end{aligned}
%\end{equation}

Salazar Celis~\cite{SalazarCelisNumericalContinued2024} proposed a heuristic for ordering nodes based on monitoring the approximation error rather than the weights: once $n$ nodes have been selected, the next node is chosen to be the candidate at which the current interpolant $r_n$ has the largest error. This greedy strategy was inspired by and mimics the adaptive selection of nodes in AAA. As Salazar Celis pointed out, this strategy rules out the halting breakdown created by an infinite weight, unless all the candidate data are already interpolated. Unlike monitoring weights, the error criterion also suggests when the iteration may be terminated due to reaching acceptable accuracy. Most importantly, Salazar Celis offered numerical evidence for challenging functions on which the greedy strategy showed no signs of significant instability. 

\begin{algorithm}
\caption{Discrete version of the greedy TCF algorithm.}
\label{alg:discrete}
\begin{algorithmic}[1]
\Require{$m$-vector of sample points $\zeta$, vector of corresponding values $\eta$, tolerance $\tau$}
\State{Let  $z,\, y,\, w$ be empty vectors}
\State{$n \gets 0$}
\State{$i \gets \arg\max_{j} |\eta_j|$}
\State{$\epsilon \gets |y_i|$}
\While{$\epsilon > \tau$ and $n < m$}
    \State{Append $\zeta_i$ to $z$ and $\eta_i$ to $y$}
    \State{$n \gets n+1$}
    \State{Calculate weight $w_n$ using~\eqref{eq:thiele-weight}}
    \State{Evaluate the interpolant $r_n$ at all unused sample points using~\eqref{eq:thiele-eval}}
    \State{$i \gets \arg\max_{j} |r_n(\zeta_j)-\eta_j|$ using~\eqref{eq:thiele-eval}} 
    \State{$\epsilon \gets |r_n(\zeta_i)-\eta_i|$}
\EndWhile
\State \Return{$z,\, y,\, w$}
\end{algorithmic}
\end{algorithm}

A more formal description of greedy TCF is given in Algorithm~\ref{alg:discrete}. Given $m$ samples of a function and an error tolerance criterion, the algorithm often halts after selecting $n \ll m$ of the sample points as nodes. The most computationally intensive step is in line 9, which requires $\O(n(m-n))$ operations to evaluate the current interpolant at all remaining candidate points. Thus, the total cost of the algorithm is $\O(m n^2)$ operations. This compares favorably with the $\O(m n^3)$ cost of AAA, which dominated by the cost of an SVD at each iteration. In making the comparison, however, we point out that the approximation accuracy is limited mainly by the rational degree, which is $\approx n/2$ for TCF and $\approx n$ for AAA, ameliorating the cost advantage of TCF. 

\section{Adaptive refinement for a continuum domain}
\label{sec:continuum}

In order to apply AAA or Algorithm~\ref{alg:discrete} to approximate a function on a domain in the complex plane, the typical first step is to discretize the domain boundary. Since the computational effort of both AAA and TCF algorithms scales just linearly with the number of sample pairs provided, and because rational approximations themselves do not require \textit{a priori} special node distributions, there may often be no need to choose the discretization too carefully. However, for functions with singularities very close to the boundary, accuracy improves greatly with nodes that cluster near the singularity, and these are the types of functions for which rational approximation offers the greatest advantage. 

\begin{figure}
    \centering
    \includegraphics[width=0.8\textwidth]{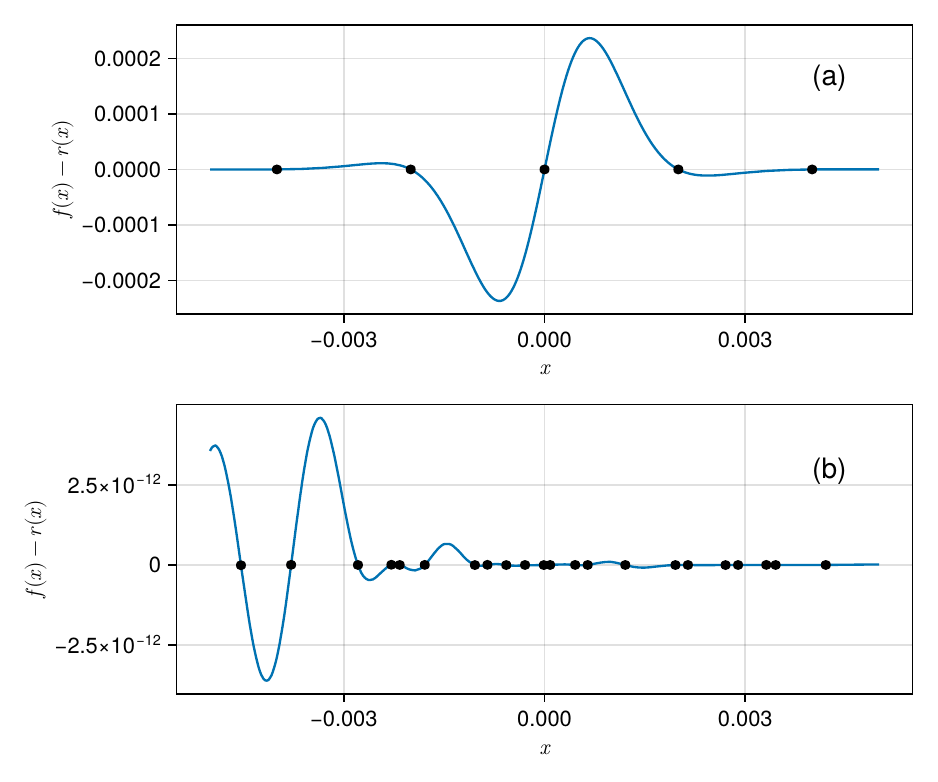}
    \caption{Error in the interpolation of $f(x)=\arctan(500x)$ on $[-1,1]$; interpolation nodes within the displayed window are shown as black points. (a) Using 1001 evenly spaced sample points in the interval $[-1, 1]$. (b) Using automatic adaptive refinement of nodes.}
    \label{fig:underresolved}
\end{figure}

\autoref{fig:underresolved} illustrates the value of adaptive refinement for the function $f(x)=\arctan(500x)$. Using 1001 evenly spaced sample points in the interval $[-1, 1]$, Algorithm~\ref{alg:discrete} finds a type $(53,53)$ interpolant with error less than $1.6\times10{-14}$ on the sample points. However, the discretization is too coarse near the step at the origin, and the error is of size $2\times 10^{-4}$ between nodes there. Algorithm~\ref{alg:continuum} below, which adaptively refines the sample set, finds a type $(48, 48)$ interpolant that is 8 orders of magnitude more accurate in the same region. 

In the extreme case of a function with a singularity on the domain, it is known that a tapered exponential distribution of nodes is ideal~\cite{TrefethenExponentialNode2021}. However, an optimal distribution is less clear for singularities that are merely nearby, and foreknowledge of the singularity structure is not always available. Hence, a method that adapts to geometry and the target function is desirable.

\begin{algorithm}
\caption{Continuum version of the greedy TCF algorithm.}
\label{alg:continuum}
\begin{algorithmic}[1]
\Require{Function $f:\mathbf{C}\to \mathbf{C}$, curve parameterization $\gamma:[0,1]\to \mathbf{C}$, tolerance $\tau$, $\nu \in \mathbf{N}$}
\State{Let $t\, z,\, y,\, w$ be empty vectors}
\State{$\epsilon \gets \infty$}
\State{$n \gets 1$}
\State{$s_* = 0$}
\While{$\epsilon > \tau$}
\State{Insert $s_*$ into $t$}
\State{$z_n \gets \gamma(s_*)$}
\State{$y_n \gets f(z_n)$}
\State{Calculate weight $w_n$ using~\eqref{eq:thiele-weight}}
\State{Evaluate vector $e_k = |r_n(\gamma(s_k)) - f(\gamma(s_k))|$ at all  $s_k \in \cup_j S_j$}
\State{$k \gets \arg \displaystyle \max_{\kappa} e_\kappa$}
\State{$s_* = s_k$}
\State{$\epsilon \gets e_k$}
\State{$n \gets n+1$}
\EndWhile
\State \Return{$t, z, y, w$}
\end{algorithmic}
\end{algorithm}

Driscoll et al.~\cite{DriscollAAARational2024} described an adaptive refinement strategy for AAA that expands the sample point set after each iteration, and the same strategy can be applied to the greedy TCF algorithm. Suppose the boundary of the domain is given by the parameterization $\gamma:[0,1]\to \mathbf{C}$, and let $\nu$ be a small positive integer. Algorithm~\ref{alg:continuum} describes an adaptive refinement strategy to produce a node set $z_k=\gamma(t_k)$, $k=1,\ldots,n$ for interpolating $f$. The algorithm maintains a sorted vector $t$ of parameter values for the selected nodes and sample sets $S_1,\dots,S_n$ interlacing them via
\begin{equation}
    S_j = \{ t_{i_j} + i\,\frac{t_{i_{j+1}} - t_{i_j}}{\nu+1}: i=1,2,\ldots,\nu \}, 
\end{equation}
with the convention that $t_{i_{n+1}}=1$. Although the terminology and initialization are different, Algorithm~\ref{alg:continuum} is the same as the refinement strategy used in~\cite{DriscollAAARational2024} for AAA.

Each evaluation of $r_n$ in Algorithm~\ref{alg:continuum} requires $\O(n)$ operations, and there are $\O(\nu n)$ such evaluations per iteration, so the total cost of the algorithm is $\O(\nu n^3)$ operations, where $n$ is the final number of nodes selected. In practice, as was done with AAA~\cite{DriscollAAARational2024}, we start with a moderately large value of $\nu=15$ and decrease it with each iteration until reaching a steady value of $\nu=3$. This allows for exponential refinement of nodes near singularities without excessive cost.

The $\O(\nu n^3)$ cost of the continuum method is to be compared to $\O(mn^3-n^3)$ for the discrete method on $m$ candidates and $\O(\nu n^4)$ for continuum AAA, with the caveat again that AAA needs only half as large a value of $n$ for the same approximation accuracy.

\section{One-division method for evaluation, derivative, and residues}
\label{sec:newalgs}

There exist well-known alternatives to~\eqref{eq:thiele-eval} for evaluation of a continued fraction. Namely, the continued fraction
\begin{equation}
    \label{eq:cfrac}
    \frac{p_n}{q_n} = b_0 + \K_{k=1}^{n} \frac{a_k}{b_k}
\end{equation}
satisfies the well-known Euler--Wallis recurrences 
\begin{equation}
    \begin{aligned}
        p_{k} &= b_{k} p_{k-1} + a_{k} p_{k-2}, & p_0 &= 1, \, p_1 = b_0, \\
        q_{k} &= b_{k} q_{k-1} + a_{k} q_{k-2}, & q_0 &= 0,\,  q_1 = 1.
    \end{aligned}
    \label{eq:recurrence}
\end{equation}
Viewed as an incremental iteration, in which $p_k$ and $q_k$ are updated from recent values given new continued fraction coefficients $a_k$ and $b_k$, this process is known to be a numerically unstable way to evaluate a continued fraction. Jones and Thron~\cite{JonesNumericalStability1974} recommended instead the tail-ordered evaluation~\eqref{eq:thiele-eval}, which has to be restarted when new coefficients are added.

The iteration~\eqref{eq:recurrence} is interpretable as the product~\cite{Milne-ThompsonCalculusFinite1933}
\begin{equation}
    \label{eq:onediv}
        \begin{bmatrix}
           p_n  \\ q_n
        \end{bmatrix}
        = 
        \begin{bmatrix}
            b_0 & 1 \\ 1 & 0
        \end{bmatrix}
        \begin{bmatrix}
            b_1 & 1 \\ a_{1} & 0
        \end{bmatrix}
        \cdots
        \begin{bmatrix}
            b_{n-1} & 1 \\ a_{n-1} & 0
        \end{bmatrix}
        \begin{bmatrix}
            b_n \\ a_{n}
        \end{bmatrix}.
\end{equation}
The unstable incremental iteration applies this product from left to right. However, if we give up incrementalism, we can instead apply it from right to left, i.e., in tail-ordered fashion. Doing so requires $n-1$ scalar multiplications, $n-1$ scalar axpy operations, and only one division, whereas~\eqref{eq:thiele-eval} requires $n-1$ divisions and $n-1$ additions. Thus,~\eqref{eq:onediv} may be faster when a floating-point division takes much longer than an axpy, as is frequently the case.

\begin{table}
    \centering
    \caption{Julia timing results in nanoseconds per evaluation of a continued fraction~\eqref{eq:cfrac} with normally distributed random coefficients $a_k,b_k$. Experiments were done in Julia 1.12.7 on an Apple M4 Pro processor.}
    \label{tab:timing-julia}
    \begin{tabular}{r|c|c|r|c|c|r}
        \multicolumn{1}{r|}{} & 
        \multicolumn{3}{c|}{Real coefficients} & \multicolumn{3}{c}{Complex coefficients} \\
        {$n$} & {Using~\eqref{eq:thiele-eval}} & {Using \eqref{eq:onediv}} & {Ratio} 
        & {Using~\eqref{eq:thiele-eval}} & {Using \eqref{eq:onediv}} & {Ratio} \\
        \hline
        % run on 2026=08-14
        25 & 82.6 & 52.5 & 1.6 & 470.4 & 76.7 & 6.1 \\
        30 & 98.9 & 63.5 & 1.6 & 568.0 & 90.6 & 6.3 \\
        35 & 118.0 & 73.4 & 1.6 & 658.0 & 101.0 & 6.5 \\
        40 & 133.9 & 83.1 & 1.6 & 754.8 & 114.6 & 6.6 \\
        45 & 147.9 & 92.4 & 1.6 & 853.1 & 126.5 & 6.7 \\
        50 & 164.3 & 101.8 & 1.6 & 952.8 & 140.9 & 6.8 \\
        55 & 178.2 & 111.5 & 1.6 & 1041 & 149.1 & 7.0 \\
        60 & 193.5 & 120.9 & 1.6 & 1137 & 162.0 & 7.0 \\
    \end{tabular}
\end{table}

\begin{table}
    \centering
    \caption{MATLAB timing results in nanoseconds per evaluation of a continued fraction~\eqref{eq:cfrac} with normally distributed random coefficients $a_k,b_k$. Experiments were done in MATLAB 26.1 on an Apple M4 Pro processor.}
    \label{tab:timing-matlab}
    \begin{tabular}{r|c|c|r|c|c|r}
        \multicolumn{1}{r|}{} & 
        \multicolumn{3}{c|}{Real coefficients} & \multicolumn{3}{c}{Complex coefficients} \\
        {$n$} & {Using~\eqref{eq:thiele-eval}} & {Using \eqref{eq:onediv}} & {Ratio} 
        & {Using~\eqref{eq:thiele-eval}} & {Using \eqref{eq:onediv}} & {Ratio} \\
        \hline
        % run on 2026=08-14
    25 & 1381 & 1348 & 1.0 & 1670 & 1482 & 1.1 \\
    30 & 1424 & 2040 & 0.7 & 1761 & 1520 & 1.2 \\
    35 & 2085 & 193.6 & 10.8 & 1924 & 380.4 & 5.1 \\
    40 & 240.1 & 203.1 & 1.2 & 719.6 & 415.8 & 1.7 \\
    45 & 254.6 & 212.2 & 1.2 & 787.8 & 454.0 & 1.7 \\
    50 & 272.0 & 227.1 & 1.2 & 868.4 & 488.1 & 1.8 \\
    55 & 288.0 & 234.5 & 1.2 & 952.6 & 527.8 & 1.8 \\
    60 & 303.7 & 249.4 & 1.2 & 1015 & 561.5 & 1.8 \\
\end{tabular}       
\end{table}

For example, \autoref{tab:timing-julia} and \autoref{tab:timing-matlab} show timing results for evaluating~\eqref{eq:cfrac} on random coefficients, in Julia and MATLAB, respectively, using both~\eqref{eq:thiele-eval} and~\eqref{eq:onediv}. In Julia, the one-division formula provides a speedup factor of about 1.6 uniformly across $25 \le n \le 60$ for real coefficients, and a speedup factor of 6--7 when the coefficients are complex. In MATLAB, those speedup factors are about 1.2 and 1.7--1.8, respectively, for $n \ge 40$.

\begin{algorithm}[thbp]
\caption{Evaluation of $p_n$, $q_n$, and their derivatives from~\eqref{eq:onediv}.}
\label{alg:onediv}
\begin{algorithmic}[1]
\State{Given: $n$-vector of nodes $z$, $n$-vector of weights $w$, evaluation point $\zeta$}
\State{$p\gets w_n$, $q\gets \zeta - z_{n-1}$}
\State{$p'\gets 0$, $q'\gets 1$}
\For{$k=n-2$ to $1$}
\State{$[p', q']\gets [w_{k+1} p' + q', p + (\zeta-z_k) p']$} (simultaneously)
\State{$[p, q]\gets [w_{k+1} p + q, (\zeta -z_k) p]$} (simultaneously)
\EndFor
\State{$[p', q']\gets [w_{1} p' + q', p']$ (simultaneously)}
\State{$[p, q]\gets [w_{1} p + q, p]$ (simultaneously)}
\State \Return{$p, q, p', q'$}
\end{algorithmic}
\end{algorithm}

The evaluation of $r_n$ from~\eqref{eq:thiele-eval} is equivalent to~\eqref{eq:onediv} with $b_k=w_{k+1}$, $a_k=z-z_k$. Algorithm~\ref{alg:onediv} summarizes the resulting iteration, with additional information returned for use in the next section. The evaluation of weights from~\eqref{eq:thiele-weight} for the case $k=n$ is also equivalent to~\eqref{eq:onediv} with
\begin{equation}
    a_k = \begin{cases}
        z_n - z_{n-k}, & 1 \le k \le n-1, \\
        y_n , & k=n,
    \end{cases} \qquad
    b_k = \begin{cases}
        0, & k=0, \\
        -w_{n-k}, & 1 \le k \le n-1, \\
        1, & k=n.
    \end{cases}
    \label{eq:weights-onediv}
\end{equation}

Computation of $r_n$ at a point via \eqref{eq:onediv} makes it easy to compute $r_n'$ as well, either by automatic differentiation software, or explicitly, as indicated in Algorithm~\ref{alg:onediv}, whose return values provide $r_n=p_n/q_n$ and $r_n'=(p_n'q_n - p_nq_n')/{q_n^2}$. The numerical stability of this formula is unknown.

Salazar~\cite{SalazarCelisNumericalContinued2024} described generalized eigenvalue algorithms for finding the roots and poles of $r_n$. Computing residues by numerical integration is not entirely straightforward when poles are crowded close together. To our knowledge, there are no published algorithms for computing residues at the poles specialized to the TCF representation. 

If $r_n=p_n/q_n$, $q_n(\zeta)=0$, $p_n(\zeta)\ne 0$, and $q_n'(\zeta)\ne 0$, then $\zeta$ is a simple pole of $r_n$ with residue
\begin{equation}
    \label{eq:residue}
    \text{Res}(r_n,\zeta) = \frac{p_n(\zeta)}{q_n'(\zeta)}.
\end{equation}
The values of $q_n'(\zeta)$ can be computed as shown in Algorithm~\ref{alg:onediv}.

\section{Computational experiments}
\label{sec:experiments}

Discrete and continuum versions of the AAA and greedy TCF algorithms have been implemented in the \textsf{RationalFunctionApproximation} package~\cite{DriscollRationalFunctionApproximationjlRational2023,DriscollRationalFunctionApproximationjl2026} in Julia.\footnote{The convergence of the Julia AAA implementation is close, but not identical, to that of the MATLAB implementation in~\cite{DriscollAAARational2024}, due to a different initialization of the node set and adjustable stopping criteria.} Here we present some experiments comparing their convergence rates and computational costs for several challenging functions defined on the interval $[-1,1]$ and on the unit circle. The experiments in this section were done using version 0.3.8 of the \textsf{RationalFunctionApproximation} package in Julia~1.12.7 on an Apple M4 Pro processor. Early termination of the iteration due to apparent stagnation, which is enabled by default in the software, was disabled for these tests. 

The convergence rates presented below---observed error as a function of rational degree---should be universal, up to fluctuations due to floating point implementation details. However, the reported timings do depend on computing architecture and environment as well as software implementation choices not fully specified by the algorithmic descriptions. This is particularly true for greedy TCF, in which the iterations are not so dominated by a single SVD operation as is the case in AAA. Nevertheless, while the exact timing values vary, we have observed consistent qualitative differences across several software and hardware versions. 

\subsection{Interval}
We considered six functions on the interval $[-1,1]$:
\begin{itemize}
    \item Singularity on the interval: $f_1(x) = \sqrt{1+x}$, $f_2(x) = |x|$
    \item Singularity near $x=0$: $f_3(x) = |x + 10^{-6}i|$, $f_5(x) = \arctan(10^6 x)$
    \item Singularity near $x=-1$: $f_4(x) = \log(x + 1 + 10^{-6})$
    \item Oscillatory: $f_6(x) = \cos(100x)$
\end{itemize}
Each was approximated using continuum variants of greedy TCF and AAA approximation up to degree 120 or a tolerance of $100$ times machine epsilon. The maximum error was measured on the union of several point sets. Define
\begin{equation}
    \label{eq:validation-sets}
    \begin{split}
        T_1 &= \{ -1 + \frac{2k}{10000}: k=0,1,\ldots,10000 \}, \\
        T_2 &= \{ 2^{-0.1k}: k=10,11,\ldots,1000 \}.
    \end{split}
\end{equation}
In order to check near all possible singularities, the validation set was defined as $T_1 \cup T_2 \cup -T_2 \cup (T_2 - 1)$. Note that not only are these points more finely distributed throughout the interval than the points used within the algorithms, but they get much closer to $x=0$ and $x=-1$ than is possible for the continuum algorithms to sample by repeated refinement in double precision. Convergence was measured as a function of the denominator degree, which allows two cases at each degree for TCF (numerator degree equal to or one greater than the denominator degree) and one case at each degree for AAA (balanced degrees only). 

\begin{figure}
    \centering
    \includegraphics[width=\textwidth]{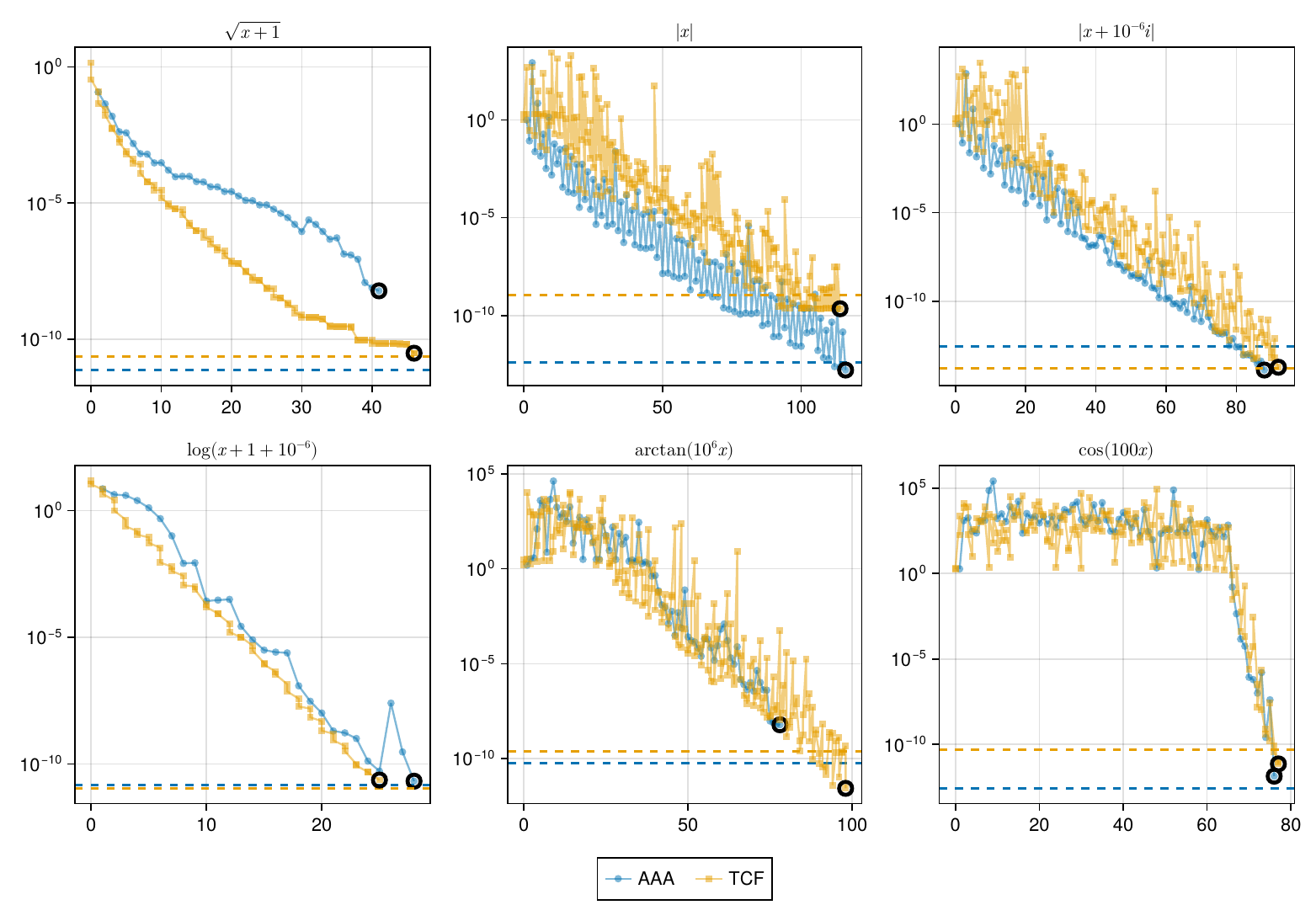}
    \caption{Comparison of TCF and AAA approximation on the interval $[-1,1]$ for six challenging functions. Each plot shows the max-norm error on a validation set described in the text as a function of the denominator degree. The black rings indicate the best approximation found for each run. The dashed lines indicate the best approximation found by discrete variants of the algorithms.}
    \label{fig:comparison-interval}
\end{figure}

\autoref{fig:comparison-interval} shows the convergence results. The dashed lines on each plot represent the error achieved by discrete variants of the algorithms using the validation set as the sample points. Thus, the discrete versions have access to exponentially refined points near any singularities, and we can check how well the continuum versions are able to find suitable nodes automatically. In all cases except $f_5$ for AAA, the continuum accuracy is comparable to or better than the discrete accuracy. Some additional experiments (not shown) indicate that this disappointing case can be improved by increasing the steady value of the refinement level $m$ in Algorithm~\ref{alg:continuum} from 3 to 4, for reasons that are unclear.

For the most part, the two methods are similar in convergence. For $f_1$, which has a singularity at an endpoint, TCF converges root-exponentially, but AAA lags behind. In this case, measurement on our finer set reveals that continuum AAA is significantly underestimating error on its adaptive test points. For $f_2$, with a singularity at $x=0$, both methods converge root-exponentially, albeit with a noticeable even/odd disparity. TCF terminates at a higher level of error in both discrete and continuum versions than AAA does. The mode of premature failure is underflow during the computation of $p_n$ and $q_n$ from~\eqref{eq:onediv} for the weight associated with adding a new node. While the underflow can easily be detected and averted by a simple rescaling, further experiments indicate that doing so offers no convergence benefit, due to an associated cumulative loss of precision during the weight iteration due to subtractive cancellation. This can partly, but not entirely, be attributed to small differences between the nodes as defined in~\eqref{eq:weights-onediv}. 

In the remaining cases, which have no singularities on the interval, both methods converge exponentially and almost identically---in the case of $f_6$, only after the function is sufficiently resolved.

\autoref{tab:comparison-interval} shows the elapsed construction times for both methods on these functions. To make the comparisons fair, each algorithm is run until reaching an error (measured on the validation set) as small as twice the greater of the best errors recorded by the two methods. This keeps the error target achievable for both methods and away from fluctuations that can occur after reaching their best results. Each construction was run 1000 times to minimize variability. The greedy TCF method was found to be faster by a factor ranging from 9 to 33. 

\begin{table}[bp]
    \caption{Timing results for finding common best interpolants from those shown in \autoref{fig:comparison-interval}.}
    \label{tab:comparison-interval}
    \centering
   \begin{tabular}{r|r|r|r|r|r}
        & \multicolumn{2}{@{}c|@{}}{\textbf{Degree}} & \multicolumn{3}{@{}c}{\textbf{Time (msec)}} \\
        {Case} & {AAA} & {TCF} & {AAA} & {TCF} & {ratio} \\
        \hline
        $\sqrt{x+1}$ & 40 & 25 & 3.3 & 0.1 & 29.9 \\
        $|x|$ & 72 & 92 & 23.5 & 2.7 & 8.8 \\
        $|x + 10^{-6}i|$ & 86 & 92 & 54.7 & 2.9 & 18.7 \\
        $\log(x + 1 + 10^{-6})$ & 28 & 25 & 1.5 & 0.1 & 11.2 \\
        $\arctan(10^6 x)$ & 75 & 78 & 34.8 & 2.0 & 17.5 \\
        $\cos(100x)$ & 76 & 77 & 31.1 & 1.5 & 20.5 \\
    \end{tabular}
\end{table}

\subsection{Unit circle} 
We considered six functions on the unit circle $|z|=1$:
\begin{itemize}
    \item Singularity at $z=-1$: $f_1(z) = \sqrt{1+z}$, $f_2(z) = |1+z|$
    \item Singularity near $z=-1$: $f_3(z) = \lvert 1 + z + 10^{-6} \rvert$, $f_4(z) = \log( 1 + z + 10^{-6} )$, $f_5(z) = \sqrt{1 + 10^6 - z^2}$
    \item Oscillatory: $f_6(z) = z^{50}$
\end{itemize}
The max error for each approximation was measured on the set 
$$
\cos(\pi T_1) \cup -\cos(\pm \pi T_2),
$$
where $T_1$ and $T_2$ are given in~\eqref{eq:validation-sets}, which allows measurement of error close to a possible singularity at $z=\pm 1$. Unlike $T_2$, however, the points in $\cos(\pi T_2)$ cannot all be distinguished in double precision, so the effective upper bound for $k$ in the definition~\eqref{eq:validation-sets} is 520 rather than 1000.

\begin{figure}
    \centering
    \includegraphics[width=\textwidth]{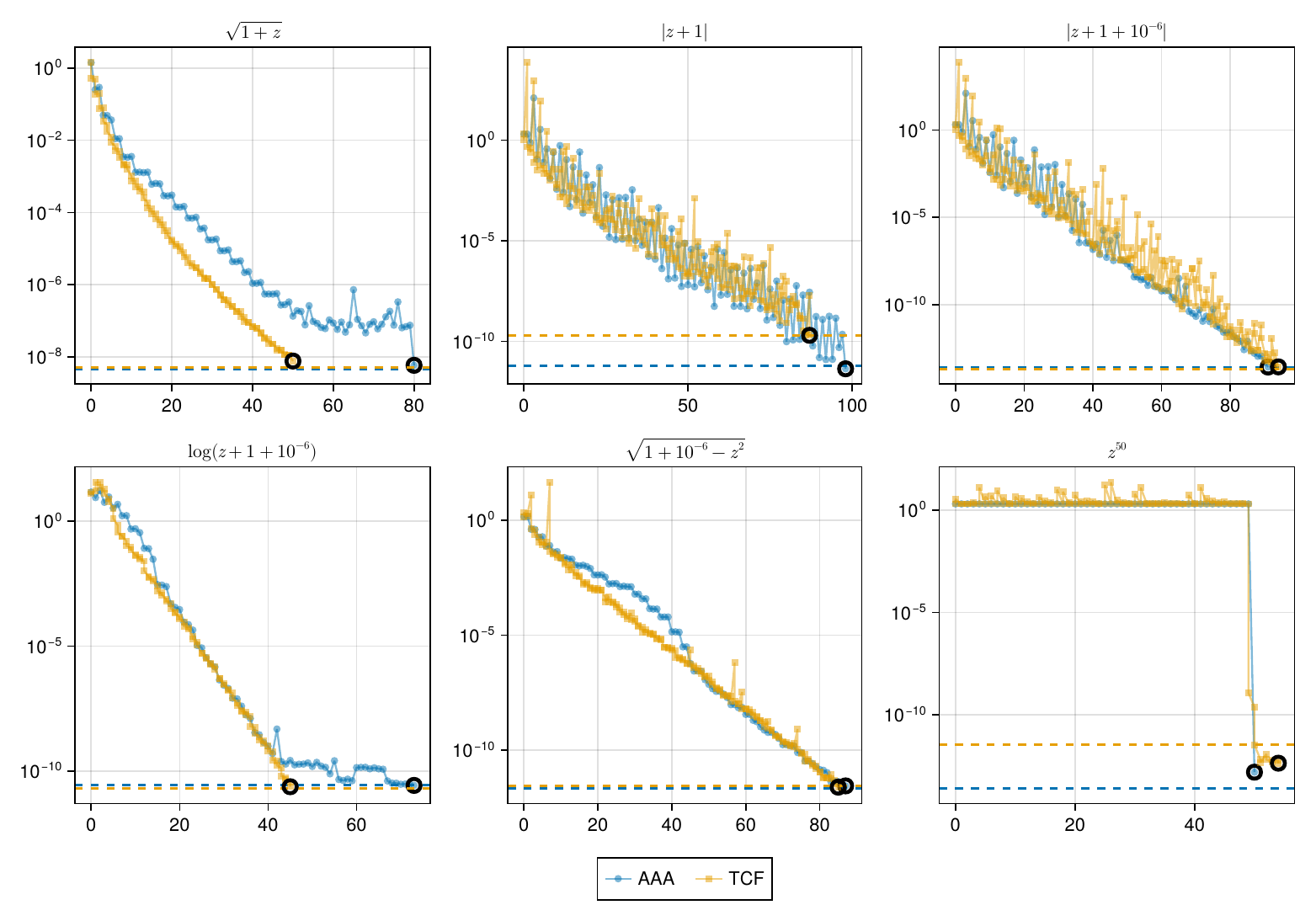}
    \caption{Comparison of greedy TCF and AAA approximation on the unit circle for six challenging functions. Each plot shows the max-norm error on a validation set described in the text as a function of the denominator degree. The black rings indicate the best approximation found for each run. The dashed lines indicate the best approximation found by discrete variants of the algorithms.}
    \label{fig:comparison-circle}
\end{figure}

The results are given in \autoref{fig:comparison-circle} and \autoref{tab:comparison-circle}. The continuum versions are able to achieve errors comparable to the discrete versions in all cases. As on the interval, the convergence rates are similar for the two methods, while the ratios of computational times range between 8 and 42. We note that the largest ratios in \autoref{tab:comparison-circle} are affected by the long, slow convergence plateaus observed for AAA in those cases; a more aggressive termination criterion would accept less accuracy for much greater speed.

\begin{table}
    \caption{Timing results for finding common best interpolants from those shown in \autoref{fig:comparison-circle}.}
    \label{tab:comparison-circle}
    \centering
    \begin{tabular}{r|r|r|r|r|r}
        & \multicolumn{2}{@{}c|@{}}{\textbf{Degree}} & \multicolumn{3}{@{}c}{\textbf{Time (msec)}} \\
        {Case} & {AAA} & {TCF} & {AAA} & {TCF} & {ratio} \\
        \hline
        $\sqrt{1 + z}$ & 80 & 46 & 60.5 & 1.4 & 42.7 \\
        $|z + 1|$ & 80 & 85 & 60.9 & 7.2 & 8.5 \\
        $|z + 1 + 10^{-6}|$ & 90 & 92 & 85.7 & 8.9 & 9.6 \\
        $\log(z + 1 + 10^{-6})$ & 56 & 44 & 19.3 & 1.3 & 15.0 \\
        $\sqrt{1 + 10^{-6} - z^2}$ & 83 & 84 & 67.5 & 6.9 & 9.8 \\
        $z^{50}$ & 50 & 51 & 14.2 & 1.7 & 8.1 \\
    \end{tabular}
\end{table}

\section{Concluding remarks}
\label{sec:conclusion}

The barycentric representation of a rational function allows specifying $n$ interpolation conditions along with $n$ weights that are used to minimize a linearized residual at sample points in a least-squares sense. The result is a rational function of type $(n-1, n-1)$. By contrast, the Thiele continued-fraction representation works only with $2n$ or $2n+1$ interpolation conditions to produce a rational approximant with denominator degree $n$. Both representations can be combined with iterative greedy node selection, either on a predetermined discrete set of candidate nodes or on a continuum domain via adaptive refinement. Even though the TCF representation only ``knows about'' the function at the selected interpolation nodes, while the barycentric representation explicitly incorporates information from additional sample points, the two methods appear to have similar global convergence rates in our tests.

Evaluation of the rational approximant and its derivative is $\O(n)$ for both representations. We have shown that by using a tail-first ordering of the classical 3-term recurrence for the TCF, we can evaluate it at any point using just one floating-point division, resulting in speed-ups of magnitude depending on the platform. Finding poles of either representation requires the solution of a generalized eigenvalue problem. We note that finding the residue at a simple pole can be done for TCF in $\O(n)$ operations through automatic differentiation. (The same should be possible for the barycentric representation, but we have not investigated this.) 

The continuum version of TCF, with no advance knowledge about function singularities, achieves results as accurate as the discrete version supplied with exponentially spaced nodes near singularities. Both versions may terminate with greater error than AAA on functions with absolute-value singularities on the domain, due to underflow/loss of significance during the TCF weight calculation. 

The TCF representation requires no linear algebra and is remarkably simple in concept and implementation. The chief appeal of greedy TCF over AAA is its lower asymptotic work requirement by a factor $\O(n)$ for construction. This advantage is offset somewhat by needing twice as many iterations to reach the same degree. In practice, our Julia implementation finds the TCF method consistently faster on all of our test functions by a factor of at least 8 and as high as 40. This speed advantage is likely to be even more pronounced in extended-precision arithmetic, when larger degrees come into play. 

In addition to the standard algorithm based on least-squares approximation, AAA offers a fast modification to find near-best approximations in the minimax sense through iterative reweighting~\cite{NakatsukasaAlgorithmReal2020}. Salazar Celis~\cite{SalazarCelisNumericalContinued2024} demonstrated a near-best method using an equilibration algorithm~\cite{HofreitherAlgorithmBest2021}, but its robustness and performance are unclear at present. AAA also has a long record of stable use for applications and successful modifications for exploiting symmetry and special problem types. The comparative simplicity and speed of TCF suggest that it is worth further investigation and development along similar lines.

\section*{Acknowledgments}
We thank Nick Trefethen and Oliver Salazar Celis for helpful discussions and suggestions.

\section{Declarations}
The authors used AI to help create code in the RationalFunctionApproximation package and for the numerical experiments. All code was reviewed manually. We assume responsibility for all the material in this work.

\bibliographystyle{siamplain}
\bibliography{Approximation}

\end{document}